\documentclass[12pt]{article}
\usepackage{amssymb}
\begin{document}
\newcommand{\nc}[2]{\newcommand{#1}{#2}}
\newcommand{\rnc}[2]{\renewcommand{#1}{#2}}
\newtheorem{proposition}{Proposition}
\newcommand{\bprop}{\medskip\begin{proposition} ~~\\ \it}
\newcommand{\eprop}{\end{proposition} \bigskip}
\nc{\ba}{\begin{array}}
\nc{\ea}{\end{array}}
\nc{\bea}{\begin{eqnarray}}
\nc{\eea}{\end{eqnarray}}
\nc{\nn}{\nonumber}
\nc{\be}{\begin{enumerate}}
\nc{\ee}{\end{enumerate}}
\nc{\beq}{\begin{equation}}
\nc{\eeq}{\end{equation}}
\nc{\bi}{\begin{itemize}}
\nc{\ei}{\end{itemize}}
\nc{\ra}{\rightarrow}
\nc{\ci}{\circ}
\nc{\lra}{\longrightarrow}
\rnc{\to}{\mapsto}
\nc{\imp}{\Rightarrow}
\rnc{\iff}{\Leftrightarrow}
\rnc{\phi}{\mbox{$\varphi$}}
\rnc{\epsilon}{\varepsilon}
\def\a{\alpha}
\def\b{\beta}
\nc{\al}{\mbox{$\alpha$}}
\nc{\hb}{\mbox{$\beta$}}
\newcommand{\Mat}{{\rm Mat}\,}
\newcommand{\half}{\frac{1}{2}}
\newcommand{\bbr}{{\bf R}}
\newcommand{\bbz}{{\bf Z}}
\newcommand{\Ci}{C_{\infty}}
\newcommand{\Cb}{C_{b}}
\newcommand{\fa}{\forall}
\newcommand{\IC}{{\mathbb C}}
\newcommand{\ID}{{\mathbb D}}
\newcommand{\IF}{{\mathbb F}}
\newcommand{\IH}{{\mathbb H}}
\newcommand{\II}{{\mathbb I}}
\newcommand{\IK}{{\mathbb K}}
\newcommand{\IM}{{\mathbb M}}
\newcommand{\IN}{{\mathbb N}}
\newcommand{\IP}{{\mathbb P}}
\newcommand{\IQ}{{\mathbb Q}}
\newcommand{\IZ}{{\mathbb Z}}
\newcommand{\IR}{{\mathbb R}}
\nc{\ot}{\otimes}
\def\otc{\otimes_{\IC}}
\def\ota{\otimes_ A}
\def\otza{\otimes_{Z(A)}}
\def\otc{\otimes_{\IC}}
%
%
\title{       
\vspace*{1cm} {\large\bf INSTANTON ALGEBRAS AND QUANTUM
$4$-SPHERES}
\vspace*{1mm}\\
\author{{\sc Ludwik D\c{a}browski}
\vspace*{-1mm}\\
\normalsize Scuola Internazionale Superiore di Studi Avanzati
\vspace*{-1mm}\\
\normalsize Via Beirut 2-4, 34014
Trieste, Italy
\vspace*{-.5mm}\\ dabrow@sissa.it
\vspace*{5mm}\\
{\sc Giovanni Landi}
\vspace*{-1mm}\\
\normalsize Dipartimento di Scienze Matematiche, Universit\`a di Trieste
\vspace*{-1mm}\\
\normalsize Via Valerio 12/b, 34127, Trieste, Italy
\vspace*{-.5mm}\\
landi@dsm.univ.trieste.it
\vspace*{.5cm}\\
}
\date{~}
}
\maketitle

\vspace{1cm}
\begin{abstract}
\noindent
We study some generalized instanton algebras which are required to
describe `instantonic complex rank $2$ bundles'. The spaces on which
the bundles
are defined are not prescribed from the beginning but rather are obtained from
some natural requirements on the instantons. They turn out to be quantum
$4$-spheres $S^4_q$, with $q\in\IC$, and the instantons are described by
self-adjoint idempotents $e$.
We shall also clarify some issues related to the vanishing of the first
Chern-Connes class $ch_1(e)$ and on the use of
the second Chern-Connes class
$ch_2(e)$ as a volume form.
\end{abstract}

\vspace*{1cm}
\begin{center}
\end{center}
\begin{center}
MS classification: 58B34, 81T75
\end{center}

\vfill\eject

\section{Introduction}
Recently there has been an intense activity on noncommutative \cite{cl00} and
quantum $4$-spheres \cite{dlm00,si,bct,bg} and instanton bundles over them.

In this paper, by generalizing the methods presented in \cite{co00,cl00}, we
search for quantum instantons. Paralleling the classical situation
\cite{aw,st,dm},
by this we just mean a complex rank
$2$  bundle, i.e. we require that the $0$th Chern-Connes class
vanishes, on some
`four dimensional space' and with not trivial characteristic classes.
Weakening the
assumptions made in
\cite{co00,cl00}, we do not require from the beginning that the
$1$th Chern-Connes class of the bundle vanishes as well.

As we shall see, the spaces on which the instanton bundle are defined are not
prescribed from the beginning but rather come {\em a posteriori}.
We could say that
first comes the bundle and then the space on which the bundle is defined.
While our procedure is completely general and could be used to produce
other quantum spaces, in this paper the resulting spaces will be quantum
$4$-sphere $S^4_q$, with $q\in\IC$. The quantum instantons will be described by
self-adjoint idempotents $e\in \Mat_4(A_q)$, with $A_q$ the
noncommutative algebra
of functions associated with (in fact, defining) the spheres $S^4_q$.
These spheres
and instantons interpolate between analogous objects recently found
in \cite{cl00}
(for $q$ such that $|q|=1$) and
\cite{dlm00} (for $q\in\IR$).

We shall also clarify some issues related to the vanishing of the $1$th
Chern-Connes class and the use of a $2$th Chern-Connes class as a
`volume form'.
It turns out that the first Chern-Connes class $ch_1(e)$ does vanish
if and only if the deformation parameter $q$ is such that $|q| = 1$.
In contrast, the second Chern-Connes class $ch_2(e)$ does not vanish
for any values
of $q$. The couple $(ch_1(e), ch_2(e))$ defines a cycle in the reduced $(b,B)$
bicomplex of cyclic homology and $ch_2(e)$ is closed, that is $b
ch_2(e) = 0$, if
and only if  $q$ is such that $|q| = 1$. It is only in the latter cases
that the class $ch_2(e)$ is `$q$-antisymmetric' and can be used as a
volume form
\cite{cl00}.

In the final section we shall make some remarks on alternative
definitions of spheres.

\vspace*{1cm}
\section{The instanton projections}

Consider first the free \mbox{$^*$-algebra} with unity
~$F = \IC [[ ~\II, \alpha_j , \beta_j , \alpha^*_j , \beta^*_j ~:~ 
j=1,2,3~ ]]$~
generated by elements $\alpha_j , \beta_j$ and their adjoints $\alpha^*_j ,
\beta^*_j$.
Then, take the following self-adjoint element $e=e^*$ in the algebra $\Mat_4(F)
\simeq \Mat_{4}(\IC ) \otimes F$,
\beq
\label{proj0} e =
\left(
\begin{array}{cc}
Q_1 & Q_2 \\ Q_{2}^{*} & Q_3
\end{array}
\right)~.
\eeq
Each of the  $Q_j$'s  is assumed to be a $2 \times 2$ matrix of
`generalized
quaternions' that is,
\beq
\label{genquat}
Q_j =
\left(
\begin{array}{cc}
\alpha_j & \beta_j \cr - q \beta^*_j & \pi \alpha^*_j
\end{array}
\right)~, ~~j = 1, 2, 3,
\eeq
and $q$ and $\pi$ are complex parameters for the time being.
Being $e$ self-adjoint also requires $Q_1= Q_1^*~, ~Q_3 = Q_3^*$~,
from which it follows that the parameter $\pi$
must, in fact, be real, that $\alpha_1 = \alpha_1^*$ and $\alpha_3 =
\alpha_3^*$,
and that
$\beta_1 = \beta_3 = 0$ (unless $q=-1$ which, for simplicity, we shall not
consider here).

\noindent
The next requirement we make is that $e$ be of rank $2$. This we implement by
requiring that its $0$th Chern-Connes class $ch_0(e)$ vanish (see later for the
definition of $ch_0$),
\beq
ch_0(e) = \left\langle \left(e - {1 \over 2} \right) \right\rangle = 0 ~.
\eeq
As a consequence we get that
\beq
(1+\pi) (\alpha_1 + \alpha_3 ) = 2 ~,
\eeq
which says that $\pi \neq -1$ and relates $\alpha_3$ with $\alpha_1$.
Summing up (and denoting $\alpha_1$ by $t$,
$\alpha_2$ by $\alpha$ and $\beta_2$ by $\beta$) up to now we have that
\bea
&&
Q_1 =
\left(
\begin{array}{cc}
t &  0 \\
0 & \pi t
\end{array}
\right)~,
~~~~~
Q_3 =
\left(
\begin{array}{cc}
{2 \over {1+\pi}} - t & 0 \\
0 & \pi( {2 \over {1+\pi}} - t)
\end{array}
\right)~,
\nn \\
&&
Q_2 = Q =
\left(
\begin{array}{cc}
\alpha & \beta \\
- q \beta^* & \pi \alpha^*
\end{array}
\right)~,
~~~~~
Q_{2}^{*} = Q^* =
\left(
\begin{array}{cc}
\alpha^* & -\bar{q} \beta \\
\beta^* & \pi \alpha
\end{array}
\right)~,
\eea
with $t=t^*$ (remember that $1+\pi \neq 0$).

\noindent
Finally we require that $e$ is idempotent as well, that is $e^2 = e$.\\
One of the consequences is that
\bea\label{idem0}
\left(
\begin{array}{cc}
t &  0 \\
0 & \pi t
\end{array}
\right) ~Q + Q~ \left(
\begin{array}{cc}
{2 \over {1+\pi}} - t & 0 \\
0 & \pi( {2 \over {1+\pi}} - t)
\end{array}
\right) = Q ~.
\eea
{}From the diagonal elements it follows that
\beq
t \alpha - \alpha t = {\pi - 1 \over \pi + 1} \alpha ~, ~~~
t {\alpha}^* - {\alpha}^* t = {1- \pi \over (\pi + 1) \pi} {\alpha}^*~,
\eeq
the consistency of which requires that $t$ commutes with
$\alpha$ and ${\alpha}^*$.
Excluding the case $\alpha = 0$ we also get that $\pi = 1$.
Then, the off-diagonal elements in (\ref{idem0}) imply that $t$
commutes with $\beta$ and ${\beta}^*$ as well.

\noindent
{}From $e^2 = e$ it also follows that
\bea
&& Q Q^* ~+~
\left(
\begin{array}{cc}
t &  0 \\
0 & t
\end{array}
\right)^2 ~-~
\left(
\begin{array}{cc}
t &  0 \\
0 & t
\end{array}
\right) = 0 ~, \nn \\
&& Q^*  Q ~+~
\left(
\begin{array}{cc}
t &  0 \\
0 & t
\end{array}
\right)^2 ~-~
\left(
\begin{array}{cc}
t &  0 \\
0 & t
\end{array}
\right) = 0 ~.
\eea
These constraints require that
\bea
&& \hb\alpha = \bar{q}\alpha\hb , ~~~
\hb^*\alpha = q \alpha\hb^* , \label{relpre} \\
&& \alpha^* \alpha + |q|^2\hb^*\hb + t^2 -t = 0 ~,
~~~~~~\alpha\alpha^* + \hb\hb^* + t^2 -t = 0 ~, \label{relpre1}\\
&& \alpha^* \alpha + |q|^2\hb\hb^* + t^2 -t = 0 ~,
~~~~~~\alpha\alpha^* + \hb^*\hb + t^2 -t = 0 ~. \label{relpre2}
\eea
Equations (\ref{relpre1}) and (\ref{relpre2}) in turn imply that
\beq
(|q|^2 - 1) (\hb\hb^* - \hb^*\hb) = 0 ~.
\eeq
Then, if $|q|^2 \neq 1$, it follows that $\hb$ must commute with
$\hb^*$. On the
other hand, if $|q|^2 = 1$, from equations (\ref{relpre1}) and
(\ref{relpre2}) it
follows directly that $\hb\hb^* = \hb^*\hb$ (and also $\alpha^* \alpha =
\alpha\alpha^*$ in this case).

\vspace*{1cm}
\section{The algebra of the quantum spheres $S^4_q$}

By slightly changing notations again, i.e. denote $t={1\over 
2}(\II-z)$ and replace $\alpha \ra
{1\over 2} \alpha , \beta \ra {1\over 2} \beta $, the construction of 
the previous section amounts
to the following.  With $q\in \IC
\setminus
\{0\}$, we consider first the free \mbox{$^*$-algebra} with unity
~$F_q = \IC [[ ~\II, \alpha , \hb , z, \alpha^* , \hb^* , z^*~ ]]$~ 
generated by three elements
$\alpha , \hb$ and $z$ (and their adjoints $\alpha^* , \hb^* , z^*$).
Then, we take the following element $e$ in the algebra $\Mat_4(F_q)
\simeq \Mat_{4}(\IC ) \otimes F_q$~,
\beq
\label{proj}
e = \half\pmatrix{ \II\! +\! z, ~~0,~~ ~~~\alpha ,~ ~~\hb~ \cr
~~~0,~~ \II\! +\! z, ~-\! q \hb^* , ~\alpha^*\cr ~~~\alpha^*, ~-\! \bar{q}\hb ,
~~\II\! -\! z, ~~0~~ \cr ~\hb^*,~ ~~~\alpha,~~ ~~~0,~~ \II\! -\! z}\ .
\eeq
By construction $e$ is self-adjoint, that is $e = e^*$.
The algebra $A_q$ of the quantum 4-sphere
$S^4_q$ is thus the quotient of
the free \mbox{$^*$-algebra} $F_q$ which results by requiring that
$e$ is idempotent as well, that is $e^2 = e$. This is
equivalent to the requirement that the generators satisfy the relations
\bea\label{s4rel}
&& z=z^*, ~~~ z\alpha = \alpha z,~~~ z\hb = \hb z, ~~~ \\
&& \hb\alpha = \bar{q}\alpha\hb , ~~~
\hb^*\alpha = q \alpha\hb^* , ~~~ \hb\hb^* = \hb^*\hb ,
\nonumber\\
&& \alpha^* \alpha + |q|^2\hb^*\hb +z^2 = \II,
~~~~~~\alpha\alpha^* + \hb\hb^* + z^2 = \nonumber \II ~.
\eea
Thus the algebra $A_q$ we are looking for is the unital \mbox{$^*$-algebra} 
generated by
the elements $\alpha , \hb$ and $z$ satisfying the relations
(\ref{s4rel}).  Here the deformation parameter
$q$ could be restricted so that $|q|  \in (0,1]$; for $q$, such that
$|q| > 1$, the transformation $q \to 1/q$, $\alpha \to \al^*$, $\hb \to
-q\hb$ and $z \to z$ yields an isomorphic sphere. \\
By restricting to $q=\exp(2\pi i\theta)$ we get the sphere $S^4_\theta$
introduced in
\cite{cl00}; while for $q\in \IR$ the present $S^4_q$ is the  same as the
one introduced in \cite{dlm00}.

When $q = 1$ the algebra of the sphere $S^4_{q}$ is commutative and
coincides with the algebra of continuous functions on the $4$
dimensional sphere $S^4$.
Thus $S^4_q$ provides a deformation of the classical sphere $S^4$.

The algebra $A_q$ can be made into a $C^*$-algebra in the usual way.
For $a\in F_q$ one defines $\| a \|$ as the supremum,
over all representations $\pi$ of $F_q$ in $B(H)$
that are {\em admissible}, in the sense that the operators
$\pi(\alpha ), \pi(\hb ), \pi(z)$ satisfy the relations (\ref{s4rel}),
of the operator norms $\| \pi (a)\|$.
Then ${\cal J} := \{ a \in F_q ~:~|a|=0 \}$
is a two-sided ideal and one obtains
a $C^*$-norm on $F_q / \cal J$.
The completion of this quotient algebra
defines a $C^*$-algebra,
which we shall denote by the same symbol $A_q$.

By using relations (\ref{s4rel}) it can be seen that the elements
$a_{kmn\ell}$, with $k\in\IZ$ and $m, n, \ell$ non negative
integers, of the form
\beq\label{basis}
a_{kmn\ell} = \cases{
\alpha^{*k} \beta^{*m}\beta^n z^{\ell}
~~~~~~\mbox{for}~k=0, 1, 2 ~\dots  \cr
\alpha^{-k} \beta^{*m}\beta^n z^{\ell}
~~~~~~\mbox{for}~k= -1, -2 ~\dots ~. }
\eeq
provides a linear basis for $A_q$.

We note that for the generic situation
when $0<|q|<1$ any character $\chi$, besides
\beq\label{s4relclq}
\chi(\al^*)=\overline{\chi(\al)}~, ~~~ \chi(\hb^*)=\overline{\chi(\hb)}~,
~~~ \chi(z^*)=\chi(z) ~,
\eeq
has to satisfy the equations
\beq
\chi(\hb) = 0 ~~{\rm and }~~ |\chi(\al)|^2+(\chi(z))^2 = 1 ~. \nn
\eeq
Thus the space of all (nonzero) characters,
which can be thought of as the space of `classical points' of $S^4_q$,
is homeomorphic to the $2$-dimensional sphere $S^2$.

Next, we describe infinite dimensional irreducible representations
of the algebra $A_q$ (for $0<|q|<1$) in $B(H)$, the algebra of 
bounded operators on a Hilbert
space $H$.  Let  $\{ \psi_n ~, ~ n = 0, 1, 2, \cdots ~\}$
be an orthonormal basis for the Hilbert space $H$.
With $\zeta \in \IC~, ~|\zeta| \leq 1$, we get two
families of representations
$\pi_{\zeta, \pm} : A_q \rightarrow B(H)$ given by
\beq\label{reps}
\begin{array}{ll}
\pi_{\zeta, \pm}(z) \, \psi_n = \pi_{\zeta, \pm}(z^*) \, \psi_n = \pm
\sqrt{ 1 - |\zeta|^2 } \, \psi_n ~,
& ~ \\
\pi_{\zeta, \pm}(\alpha)
\, \psi_n = \zeta \, \sqrt{ 1 - |q|^{2(n+1)} } \,
\psi_{n+1} ~,
&
\pi_{\zeta, \pm}(\alpha^*) \, \psi_n = \bar{\zeta} \,
\sqrt{ 1 -
|q|^{2n} }
\, \psi_{n-1} ~, \\
\pi_{\zeta, \pm}(\beta) \, \psi_n =
\zeta \, \bar{q}^n \, \psi_n ~,
& \pi_{\zeta, \pm}(\beta^*) \, \psi_n
= \bar{\zeta} \, q^n \, \psi_n ~,
\end{array}
\eeq
In fact, for $\zeta$ such that $|\zeta| = 1$, the two
representations $\pi_{\zeta,+}$ and $\pi_{\zeta,-}$ are identical so
that the representations are parametrized by points on
a classical sphere $S^2$, similarly to what happens for one dimensional
representations (characters) as described before.

\vspace*{1cm}
\section{Chern-Connes classes}
The self-adjoint idempotent $e$ given by (\ref{proj}) is clearly
an element in the matrix algebra 
$\Mat_4(A_q) \simeq \Mat_{4}(\IC ) \otimes A_q$.
It naturally acts on the
right free $A_q$-module $A_q^4 = A_q\otimes \IC^4$, and one gets as its
range a projective module of finite type which may be thought of as
the module of `sections of a vector bundle over $S^4_q$ '. The module
$e A_q^4$ is a deformation of the classical instanton bundle over
$S^4$: for $q=1$, the module $e A_q^4$ is the module of sections of the
complex rank two, instanton bundle over $S^4$ \cite{at}. 

We compute now the Chern-Connes Character 
of the module $e A_q^4$ for a generic value of the deformation
parameter $q$.
If $\langle \ \rangle$ is
the projection on the commutant of $4 \times 4$ matrices (a partial trace),
up to normalization the component of the (reduced) Chern-Connes Character
are given by the formul{\ae}
\beq
ch_n(e) = \left\langle \left(e - {1 \over 2} \right) \underbrace{e
\ot \cdots
\ot e}_{2n}
\right\rangle ~,
\eeq
and they are elements of the tensor product
\beq
A_q \ot \underbrace{\bar{A_q} \ot \cdots \ot \bar{A_q}}_{2n} ~,
\eeq
where $\bar{A_q} = A_q / \IC \II$ is the quotient of the algebra $A_q$ by the
scalar multiples of the unit $\II$. \\
The crucial property of the components ${\rm ch}_n (e)$ is that they 
define a {\it cycle} in the
$(b,B)$ bicomplex of cyclic homology \cite{co85,co94,L,co00a},
\begin{equation}
B \, {\rm ch}_n (e) = b \, {\rm ch}_{n+1} (e) \, . \label{Bb}
\end{equation}
The operator $b$ is defined by
\beq\label{bi}
b (a_0 \ot a_1 \ot \cdots \ot a_m ) =
\sum_{j=0}^{m-1} (-1)^j a_0 \ot \cdots \ot a_j a_{j+1} \ot \cdots \ot a_m
\, + (-1)^{m} a_m a_0 \ot a_1 \ot \cdots \ot a_{m-1} \,
\eeq
while the operator $B$ is written as
\beq\label{Bi}
B =  B_0 A\, ,
\eeq
where
\bea\label{Bi0}
&& B_0 (a_0 \ot
a_1 \ot
\cdots \ot a_m ) = \II \ot a_0 \ot a_1 \ot \cdots
\ot a_m
\\
&& A
(a_0 \ot a_1 \ot \cdots \ot a_m ) = { 1 \over m + 1} \sum_{j=0}^m 
(-1)^{mj} a_j \ot a_{j+1} \ot \cdots \ot a_{j-1} \, , \label{Ai}
\eea
with the obvious cyclic identification $m+1 = 0$.
To be precise, in formul{\ae} (\ref{bi}), (\ref{Bi0}) and (\ref{Ai}),
all elements in the tensor products but the first one should be taken modulo
complex multiples of the unit $\II$, that is one has to project onto
$\bar{A_q} = A_q / \IC \II$.

\bigskip
\noindent
The fact that $ch_0(e) = 0$ has been imposed from the very beginning 
and was one of
the conditions that lead to the projections (\ref{proj}). As already remarked,
this could be interpreted as saying that the idempotent and the corresponding
module (the `vector bundle') has complex rank equal to $2$.

\bigskip
\noindent
Next one finds,
\bea
ch_1(e) &=& \left\langle \left(e - {1 \over 2} \right) \ot e \ot e
\right\rangle
\\ &=& {1 \over 8} (1-|q|^2) \Big\{
z \ot (\beta \ot \beta^* - \beta^* \ot \beta ) \nn \\
&& ~~~~~~~~~~~~~~~~~~~~~~~~
+ \beta^* \ot (z \ot \beta - \beta \ot z) + \beta \ot (\beta^* \ot z -
z \ot \beta^*) \Big\}
   ~. \nn
\eea
It is straightforward to check that
\beq
b ch_1(e) = 0 = B ch_0(e)
\eeq
Furthermore, we have the following
\bprop
Given the projections (\ref{proj}) its first Chern-Connes class
vanishes, $ch_1(e)=0$, if and only if the deformation parameter $q$
is such that
$|q|=1$.
\eprop
This result matches the analogous one found in \cite{cl00}.

\bigskip
\noindent
Finally,
\beq
ch_2(e) = \left\langle \left(e - {1 \over 2} \right) \ot e \ot e\ot e \ot
e \right\rangle
\eeq
is the sum of five terms
\beq\label{vol}
ch_2(e) = {1\over 32} \, \Big( z \, \ot \, c_z + \a \, \ot \, c_{\a} + \a^*
\, \ot \,
c_{\a^*} +
\b \, \ot \, c_{\b} + \b^* \, \ot \, c_{\b^*} \Big) \, , \label{hocy}
\eeq
with
\begin{eqnarray}\label{gamz}
c_z &=& (1-|q|^4) \, (\b \, \ot \, \b^* \, \ot \, \b \, \ot \, \b^*
-  \b^* \, \ot \, \b \, \ot \, \b^* \, \ot \, \b) \\
&+& (1-|q|^2) \, \Big\{
z \, \ot \, z \, \ot \, (\b \, \ot \, \b^* -  \b^* \, \ot \, \b)
+ (\b \,
\ot \, z \, \ot \, z \, \ot \, \b^* -  \b^*  \, \ot \, z \, \ot \, z \, \ot
\, \b) \nn \\
&\, & \, \, \, \, \, \, \, \, \, \, \, \, \, \, \, \,\,
\, \,\, \, \,\, \, \,
+  \, (\b \, \ot \, \b^* - \b^* \, \ot \, \b)  \, \ot
\,
z \, \ot \,
z  + z \, \ot \, (\b \, \ot \, \b^* - \b^* \, \ot \, \b) \,
\ot \, z
\nn
\\ &\, & \, \, \, \, \, \, \,\, \, \, \, \,
- z \, \ot \, (\b \,
\ot \, z \, \ot \, \b^* -  \b^* \, \ot \, z \, \ot \, \b)
- (\b \,
\ot \, z
\, \ot \, \b^* -  \b^* \, \ot \, z \, \ot \, \b) \, \ot \, z
\,
\Big\} \nn \\
&+& (\a \, \ot \, \a^* - |q|^2 \, \a^* \, \ot \,
\a)
\, \ot
\, (\b \, \ot \, \b^* -  \b^* \, \ot \, \b)
\nonumber \\
&\, & \, \, \, \,
\, \, \, \, \, \,
+ \, ( \b \, \ot \, \b^* -  \b^*
\, \ot \, \b)
\, \ot \,
( \a \, \ot \, \a^* - |q|^2 \, \a^* \, \ot \,
\a)
\nonumber \\
&+& \, (\b
\, \ot \, \a  - \bar{q} \, \a \, \ot \,
\b)
\, \ot \, (\a^* \, \ot \, \b^*
- q \, \b^* \, \ot \, \a^*)
\nonumber
\\
&\, & \, \, \, \, \, \, \, \, \,
\,
+ \, (\a^* \, \ot \, \b^* - q
\, \b^* \, \ot \, \a^*)
\, \ot \, (\b \,
\ot \, \a  - \bar{q} \, \a
\, \ot \, \b)
\nonumber \\
&+& \, ( \a^* \, \ot \, \b - \bar{q} \,
\b \, \ot \, \a^*)
\, \ot \, (q \,  \a \, \ot \, \b^*
\, -  \b^* \,
\ot \, \a)
\nonumber \\
&\, & \, \, \, \, \, \, \, \, \, \,
+
\, (q
\,  \a \, \ot \, \b^* \, -   \b^* \, \ot \, \a)
\, \ot \, ( \a^*
\,
\ot \, \b - \bar{q} \,  \b \, \ot \, \a^*) \, ;
\nonumber
\end{eqnarray}

\begin{eqnarray}\label{gamalp}
c_\a &=& (z
\,
\ot \, \a^* - \a^* \, \ot \, z)
\, \ot \, (\b^* \, \ot \, \b -  \b
\, \ot
\, \b^*)
\\ &\, & \, \, \, \, \, \, \, \, \, \,
+ \, |q|^2 \,
(\b^* \, \ot
\, \b -  \b \, \ot \, \b^*)
\, \ot \, (z \, \ot \, \a^*
- \a^* \, \ot \, z)
\nonumber \\
&+& \, \bar{q} \, (z \, \ot \, \b  -
\b \, \ot \, z)
\, \ot
\, (\a^* \, \ot \, \b^* - q \, \b^* \, \ot \,
\a^*)
\nonumber \\
&\, & \,
\, \, \, \, \, \, \, \, \,
+ \, (\a^* \,
\ot \, \b^* - q \, \b^* \, \ot \,
\a^*)
\, \ot \, (z \, \ot \, \b  -
\b \, \ot \, z)
\nonumber \\
&+& \, q \,
(\b^* \, \ot \, z \, - z \,
\ot \, \b^*)
\, \ot \, ( \a^* \, \ot \, \b -
\bar{q} \,  \b \, \ot \,
\a^*)
\nonumber \\
&\, & \, \, \, \, \, \, \, \,
\, \,
+ \, ( \a^* \,
\ot \, \b - \bar{q} \, \b \, \ot \, \a^*)
\, \ot \,
(\b^* \, \ot \, z
\, - z \, \ot \, \b^*) \,
;
\nonumber
\end{eqnarray}
\begin{eqnarray}\label{gamalpbar}
c_{\a^*}
&=&
|q|^2
\, (z \, \ot \, \a - \a \, \ot \, z)
\, \ot \, (\b \, \ot \, \b^*
-
\b^* \, \ot \, \b)
\\ &\, & \, \, \, \, \, \, \, \, \, \,
+ \, (\b
\,
\ot \, \b^* -  \b^* \, \ot \, \b)
\, \ot \, (z \, \ot \, \a - \a
\, \ot \, z)
\nonumber \\
&+& \, (\b^* \, \ot \, z - z \, \ot \, \b^*
)
\, \ot \, (\b \, \ot \, \a - \bar{q} \, \a \, \ot \, \b)
\nonumber
\\
&\, & \, \, \, \,
\, \, \, \, \, \,
+ \, q \, (\b \, \ot \, \a -
\bar{q} \, \a \, \ot \, \b)
\, \ot \, (\b^* \, \ot \, z - z \, \ot \,
\b^* )
\nonumber \\
&+& \, (z \, \ot \, \b  \, - \b \, \ot \, z)
\,
\ot \, (\b^* \, \ot \, \a - q \,  \a
\, \ot \, \b^*)
\nonumber \\
&\,
& \, \, \, \, \, \, \, \, \, \,
+ \,
\bar{q} \, (\b^* \, \ot \, \a -
q \,  \a \, \ot \, \b^*)
\, \ot \, (z \,
\ot \, \b  \, - \b \, \ot \,
z) \,
;
\nonumber
\end{eqnarray}

\begin{eqnarray}\label{gambet}
c_\b
&=&
(1-|q|^4)
\, \Big[ (\b^* \, \ot \, z - z \, \ot \, \b^*)  \, \ot \, \b
\,  \ot
\, \b^* +  \b^* \, \ot \, \b \, \ot (\b^* \, \ot \, z - z \, \ot
\,
\b^*) \Big] \\
&+& (1-|q|^2) \, \Big\{
\b^* \, \ot \, z \, \ot \,
z \, \ot
\, z
- z \, \ot \, \b^* \, \ot \, z \, \ot \, z  \nn \\
&\,
& \, \, \, \,
\, \, \, \, \, \, \, \, \, \, \, \, \, \, \, \,
\, \,
\, \, \, \, \, \, \,
\, \, \, \, \, \, \, \, \, \, \,
\, \, \, \, \,
\, \, \, \, \, \, \, \, \,
\, \, \, \, \, \,
+ z \, \ot \, z \, \ot
\, \b^* \, \ot \, z
- z \, \ot \,
z \, \ot \, z\ot \, \, \b^* \,
\Big\} \nn \\
&+& (\b^* \, \ot \, z - z \,
\ot \, \b^*)
\, \ot \, (\a
\, \ot \, \a^* - |q|^2 \, \a^* \, \ot \, \a)
\nn
\\ &\, & \, \, \, \,
\, \, \, \, \, \,
+ \, (\a \, \ot \, \a^* - |q|^2 \,
\a^* \, \ot \,
\a)
\, \ot \, (\b^* \, \ot \, z - z \, \ot \,
\b^*)
\nonumber \\
&+&
\, (\a \, \ot \, z - z \, \ot \, \a)
\, \ot \, (\a^*
\, \ot \, \b^* -
q \, \b^* \, \ot \, \a^*)
\nonumber \\
&\, & \, \, \, \,
\, \, \, \,
\, \,
+ \, \bar{q} \, (\a^* \, \ot \, \b^* - q \, \b^* \, \ot
\,
\a^*)
\, \ot \, (\a \, \ot \, z - z \, \ot \, \a)
\nonumber \\
&+&
\,
(\b^* \, \ot \, \a - q \, \a \, \ot \, \b^*)
\, \ot \, (\a^* \,
\ot \, z \,
- z \, \ot \, \a^*)
\nonumber \\
&\, & \, \, \, \, \, \,
\, \, \, \,
+ \,
\bar{q} \, (\a^* \, \ot\, z \, - z \, \ot \, \a^*)
\, \ot \, (\b^*
\,
\ot
\, \a - q
\, \a \,
\ot \, \b^*) \,
;
\nonumber
\end{eqnarray}

\begin{eqnarray}\label{gambetbar}
c_{\b^*}
&=&
(1-|q|^4)
\, \Big[ (z \, \ot \, \b - \b \, \ot \, z)  \, \ot \, \b^* \,
\ot
\,
\b
+  \b \, \ot \, \b^* \, \ot (z \, \ot \, \b - \b \, \ot \, z
)
\Big] \\
&+& (1-|q|^2) \,
\Big\{ - \b \, \ot \, z \, \ot \, z \,
\ot \, z
+ z \, \ot \, \b \, \ot \, z \, \ot \, z
\nonumber \\
&\, &
\, \, \, \, \, \,
\, \, \, \, \, \, \, \, \, \, \, \, \, \, \, \,
\,
\, \, \,
\, \, \, \,\,
\, \, \, \, \, \, \, \, \,
- z \, \ot \, z
\,\ot \, \b \, \ot \, z
+ z \,
\ot \, z \, \ot \, z\ot \, \, \b \,
\Big\} \nn \\
&+& (z \, \ot \, \b - \b
\, \ot \, z)
\, \ot \, (\a \,
\ot \, \a^* - |q|^2 \, \a^* \, \ot \, \a)
\nn
\\ &\, & \, \, \, \, \,
\, \, \, \, \,
+ \, (\a \, \ot \, \a^* - |q|^2 \,
\a^* \, \ot \,
\a)
\, \ot \, (z \, \ot \, \b - \b \, \ot \, z)
\nonumber
\\
&+& \, q
\, (z \, \ot \, \a^* - \a^* \, \ot \, z)
\, \ot \, (\b \, \ot
\, \a -
\bar{q} \, \a \, \ot \, \b)
\nonumber \\
&\, & \, \, \, \, \, \,
\,
\, \, \,
+ \, (\b \, \ot \, \a - \bar{q} \, \a \, \ot \, \b)
\,
\ot \, (z
\, \ot \, \a^* - \a^* \, \ot \, z)
\nonumber \\
&+& q \,
(\a^* \, \ot \, \b
- \bar{q} \, \b \, \ot \, \a^*)
\, \ot \, (z \,
\ot \, \a - \a \, \ot \, z) \,
\nonumber \\
&\, & \, \, \, \, \, \,
\, \, \, \,
+ \, \, (z \, \ot \, \a
- \a \, \ot \, z)
\, \ot \, (\a^*
\, \ot \, \b - \bar{q} \, \b \, \ot \,
\a^*) \, .
\nonumber
\end{eqnarray}

\bigskip
\noindent
By using the relations
(\ref{s4rel}) for our algebra, and remembering
that we need to
project on $\bar{A_q}$ in all terms of the tensor product
but
the
first one, a long (one needs to compute 750 terms) but
straightforward
computation gives
\bea\label{bch2}
&& b ch_2(e) = {1 \over 16} \, (1-|q|^2) \, \Big\{
\, \II \, \ot \, z \, \ot \, ( \b \, \ot \, \b^*
- \b^* \, \ot \, \b )
\\ && \, \, \, \, \, \, \, \, \, \, \, \, \, \, \, \, \, \, \, \, \, \, \,
\, \,
\,\, \, \, \, \, \, \, \, \, \, \, \, \, + \, \II \, \ot \, \b \, \ot \, (\b^*
\ot \, z \, - z \ot \b^*)
   + \, \II \,
\ot \, \b^* \ot \, (z \, \ot \, \b - \b \ot \, z ) \, \Big\} \nn
\eea
and this is exactly equal to $B ch_1(e)$. \\
Again as in \cite{cl00}, we shall have the following
\bprop
Given the projection in (\ref{proj}), the cycle $ch_2(e)$ is closed, i.e.
$bch_2(e)=0$, if and only if the deformation parameter $q$ is such
that $|q|=1$.
Then, the resulting class $ch_2(e)$ is `$q$-antisymmetric' and can be used as a
`volume form'.
\eprop

\vspace*{1cm}
\section{Final remarks}

In the present paper, we have searched for a general
algebra $A$ for which the element $e\in \Mat_4(A)$
is an hermitian rank $2$ projection of the form (\ref{proj0}).
As a result we have obtained a family of quantum 4-spheres $S^4_q$,
$q\in \IC$,
which turn out to be a suspension of a family of quantum 3-spheres $S^3_q$.
In the two special cases $q\in \IR$ and $|q| = 1$,
we obtain the two families found in \cite{dlm00} and \cite{cl00},
respectively.
Moreover, our instanton projections $e$ also specialize at the same time
to the projections presented therein.
This explains in  which sense these two special families are
related by analytic continuation of the deformation parameter.

While this work was being completed,
there has been the papers \cite{si} and its generalization \cite{bg},
where other families $S^4_{q,\theta}$, $q, \theta \in \IR$
and $S^4_{p,q,s}$, $|p|=1, q, s \in \IR$
of noncommutative $4$-spheres have been introduced.
In the appropriate limits (when  $|q|=1$)
they exactly reduce to
the quantum $4$-sphere $S^4_\theta$ of \cite{cl00}.
However, the projections defined therein do not specialize
correspondingly to those of \cite{cl00}.

Provided that we could perform the polar decomposition
$\hb = |\hb |~{\rm phase}(\hb)$ of $\hb$,
the change of variables $q^{\prime}\mapsto |q|, e^{2i\theta} \mapsto 
{\bar q}/|q|,
\alpha^{\prime} \mapsto  \alpha  , \hb^{\prime} \mapsto  |q|~|\hb|$ and
$U \mapsto  z~ {\rm phase}(\hb)$,
seems to suggest that the family  $S^4_{q,\theta}$ (and then $S^4_{p,q,s}$)
should not be too much different from our $S^4_{q}$.
This is, however, not the case and these families are not equivalent
as it is clear e.g. from the fact that
the related spaces of characters are different.

We also remark that the form of the projections $e$ in \cite{si}
(and in \cite{bg}) is different from ours.
The same holds for $e^{\prime}$ of \cite{si}
but $e^{\prime}$, being independent of ${\rm phase}(\hb)$,
is clearly even nonequivalent.
In fact $e^{\prime}$, rather than on a noncommutative
$4$-sphere, lives on a $3$-sphere and thus,
not surprisingly,
corresponds to some trivial bundle (with vanishing $ch_k(e)$, $k=0, 1, 2$).
(The same observations apply to $\tilde{e}$ in \cite{bg}).

\bigskip\bigskip

\bigskip\bigskip \bigskip\bigskip

\vfill\eject

\end{document}